\documentclass[a4paper, cleveref, autoref, thm-restate]{lipics-v2021}

\hideLIPIcs

\usepackage[T1]{fontenc}
\usepackage{amsmath,amssymb}

\usepackage{amsthm}

\usepackage{graphicx}
\usepackage{hyperref}
\newtheorem{question}{Question}
\crefname{question}{Question}{Questions}
\Crefname{question}{Question}{Questions}

\usepackage{color}

\usepackage{cite}

\newcommand{\T}{\mathcal{T}}
\newcommand{\twG}{{\rm tw}(G)}
\newcommand{\tw}{{\rm tw}}

\newcommand{\treeaG}{{\rm tree\textnormal{-}}\alpha(G)}
\newcommand{\treea}{{\rm tree}\textnormal{-}\alpha}
\newcommand{\treex}{{\rm tree}\textnormal{-}\chi}
\newcommand{\treexG}{{\rm tree\textnormal{-}}\chi(G)}

\nolinenumbers
\title{On the Relation Between Treewidth, Tree-Independence Number, and Tree-Chromatic Number of Graphs}
\author{Alex Koutsoutis}{Texas A\&M University, USA}{akoutsoutis@tamu.edu}{}{}

\author{Kilian Krause}{Karlsruhe Institute of Technology, Germany}{KilianKrause@outlook.de}{}{}

\author{Chun-Hung Liu}{Texas A\&M University, USA}{chliu@tamu.edu}{}{Partially supported by NSF under CAREER award DMS-2144042.}

\author{Mirza Redzic}{Karlsruhe Institute of Technology, Germany}{mirza.redzic}{https://orcid.org/0009-0001-7509-1686}
{supported by the Deutsche Forschungsgemeinschaft (DFG, German Research Foundation) -- 462679611.}

\author{Torsten Ueckerdt}{Karlsruhe Institute of Technology, Germany}{torsten.ueckerdt@kit.edu}{https://orcid.org/0000-0002-0645-9715}{supported by the Deutsche Forschungsgemeinschaft (DFG, German Research Foundation) -- 532213793.}

\authorrunning{A. Koutsoutis et al.}

\titlerunning{On the relation between treewidth, $\treea$ and $\treex$ of a graph}

\keywords{Tree-independence number, Tree-chromatic number, Treewidth.}

\ccsdesc[500]{Mathematics of computing~Extremal graph theory}

\begin{document}

\maketitle

\begin{abstract}
Dallard, Milani\v{c}, and \v{S}torgel (JCTB, 2024) asked whether every graph $G$ satisfies that $\twG+1 \leq \treeaG \cdot \treexG$, where $\tw(G)$, $\treea(G)$ and $\treex(G)$ are the treewidth, the tree-independence number and the tree-chromatic number of $G$, respectively.
We provide a negative answer for this question in a strong form: for every function $f\colon {\mathbb N} \rightarrow {\mathbb N}$, there exists a graph $G$ such that $\tw(G) > \treea(G) \cdot f(\treex(G))$.
On the other hand, we observe that $\tw(G)+1 \leq \treea(G)^2 \cdot \treex(G)$ for every graph $G$.
%
%    We investigate two recently introduced graph parameters, both of which measure the complexity of the tree decompositions of a given graph.
%    The treewidth $\twG$ of a graph $G$ measures the largest number of vertices required in a bag of every tree decomposition of $G$.
%    Similarly, the tree-independence number $\treeaG$ and the tree-chromatic number $\treexG$ measure the largest independence number, respectively the largest chromatic number, required in a bag of every tree decomposition of $G$.
%    
%    Recently, Dallard, Milani\v{c}, and \v{S}torgel asked (JCTB, 2024) whether for all graphs $G$ it holds that $\twG+1 \leq \treeaG \cdot \treexG$.
%    We answer that question in the negative, providing two infinite families of graphs $G$ with $\twG+1 > \treeaG \cdot \treexG$.
%    In our best lower bound $\twG$ is roughly $2\cdot \treeaG \cdot \treexG$.
%
%    For our approach we introduce $H$-subdivisions of $G$, which may be of independent interest, and give sharp bounds on the treewidth, tree-independence number, and tree-chromatic number of such graphs.
%    Additionally, we give upper bounds on $\twG$ in terms of $\treeaG$ and $\treexG$ that hold for all graphs $G$.
\end{abstract}

\section{Introduction}
\label{sec:introduction}

Tree-decompositions constitute a central pillar in algorithmic and structural graph theory.
While originally introduced by Halin~\cite{Hal-76} in 1976, this notion was rediscovered independently and developed thoroughly by Robertson and Seymour within their Graph Minors project. %Robertson and Seymour~\cite{RS-84} within their Graph Minors project.

Given a graph $G$, a \emph{tree-decomposition} of $G$ is a pair $\T = (T,\{X_t\}_{t\in V(T)})$ consisting of a tree $T$ and for each node $t$ in $T$ a set $X_t \subseteq V(G)$, called a \emph{bag} at $t$, such that: % the following hold:
\begin{itemize}
    %\item for every vertex $v \in V(G)$, there exists a node $t \in V(T)$ with $v \in X_t$
    \item for every edge $uv \in E(G)$, there exists a node $t \in V(T)$ with $\{u,v\} \subseteq X_t$, and
    \item for every vertex $v \in V(G)$, the set $\{t \in V(T): v \in X_t\}$ induces a non-empty connected subgraph of $T$. %$X^{-1}(v) = \{t \in V_T \mid v \in X_t\}$ of all nodes whose bag contains $v$ forms a connected subgraph of $T$
\end{itemize}
%Every graph $G$ admits a tree decomposition: taking a tree on a single node whose bag contains all vertices of $G$. % is called the \emph{trivial tree decomposition} $\T = (T = (\{t\},\emptyset), \{X_t = V\})$ of $G$.

Due to various applications, tree-decompositions whose bags induce ``tractable'' subgraphs are considered (e.g., \cite{BBEGLPS-24,DMS-24,DVORAK2025103631,DM-18,Distel_2025,DOURISBOURE20072008,DK-14,Liu-25,LNW-24,ROBERTSON1986309,ROBERTSON200343,Sey-16,Yol-18}).
There are different ways to measure the ``tractability'' of the subgraphs induced by bags.
One way to do so is considering certain graph parameters of those subgraphs.
For a graph parameter $p$ and a tree-decomposition $\T = (T,\{X_t\}_{t\in V(T)})$ of a graph $G$, we define $p(\T)$ to be $\max_{t \in V(T)}p(G[X_t])$; we define $${\rm tree}\textnormal{-}p(G)=\min_{\T'}p(\T'),$$ where the minimum is over all tree-decompositions $\T'$ of $G$.

The most classical example is the treewidth of a graph $G$.
%Formally, the \emph{width} of a given tree-decomposition is the maximum size of a bag minus 1.
The \emph{width} of a graph is its number of vertices minus 1.
Hence the width of a given tree-decomposition is the maximum size of a bag minus 1, and the tree-width of a graph $G$, denoted by $\tw(G)$, is the minimum width of a tree-decomposition of $G$.
Tree-width is also written as \emph{treewidth} in the literature.
%The \emph{treewidth} of a graph $G$, denoted by $\tw(G)$, is the minimum width of a tree-decomposition of $G$.
Graphs with bounded treewidth have numerous nice properties.
For example, treewidth determines the minimum size of balanced separators \cite{ROBERTSON1986309}, and every property that can be expressed in monadic second order logic can be decided in linear time for graphs of bounded treewidth \cite{COURCELLE199012}.

Another example is the tree-independence number, which was introduced independently by Yolov \cite{Yol-18} (under the name $\alpha$-treewidth) and Dallard, Milani\v{c}, and \v{S}torgel~\cite{DMS-24}, measuring subgraphs induced by bags by their independence number.
The \emph{independence number} of a graph $H$, denoted by $\alpha(H)$, is the maximum size of a stable set (i.e., a set of pairwise non-adjacent vertices).
%For a tree-decomposition $\T=(T,\{X_t\}_{t \in V(T)})$ of a graph $G$, we define $\alpha(\T)$ to be $\max_{t \in V(T)}\alpha(G[X_t])$. 
%The \emph{tree-independence number} of a graph $G$, denoted by $\treea(G)$, is $\min_\T\alpha(\T)$ over all tree-decompositions $\T$ of $G$.
Hence the tree-independence number of a graph $G$ is the minimum $k$ such that there exists a tree-decomposition of $G$ such that every subgraph induced by a bag has independence number at most $k$.
Tree-independence number enjoys several features of treewidth. %such as the existence of certain kind of balanced separators \cite{CHUDNOVSKY202674} and the 
For example, tree-independence number is effective for solving algorithmic problems related to independence number \cite{DMS-24,lima_et_al:LIPIcs.ESA.2024.85}.

The third example considered in this paper is the tree-chromatic number, introduced by Seymour \cite{Sey-16}, measuring subgraphs induced by bags by their chromatic number.
The \emph{chromatic number} of a graph $H$, denoted by $\chi(H)$, is the smallest integer $k$ such that $V(H)$ can be partitioned into $k$ stable sets.
%For a tree-decomposition $\T=(T,\{X_t\}_{t \in V(T)})$ of a graph $G$, we define $\chi(\T)$ to be $\max_{t \in V(T)}\chi(G[X_t])$. 
%The \emph{tree-chromatic number} of a graph $G$, denoted by $\treex(G)$, is $\min_\T\chi(\T)$ over all tree-decompositions $\T$ of $G$.
Hence the tree-chromatic number of a graph $G$ is the minimum $k$ such that there exists a tree-decomposition of $G$ such that every subgraph induced by a bag has chromatic number at most $k$.
Unlike treewidth and tree-independence number, tree-chromatic number is more subtle and less understood.

Inspired by the simple but very useful inequality $|V(G)| \leq \alpha(G) \cdot \chi(G)$, which holds for every graph $G$, Dallard et al.\ \cite{DMS-24} proposed the following question addressing its tree-decomposition analog.

\begin{question}[{Dallard, Milani\v{c}, and \v{S}torgel~\cite[Question 8.4]{DMS-24}}]{\ \\}\label{que:main-question}%
    Does every graph $G$ satisfy $\twG+1 \leq \treeaG \cdot \treexG$?
\end{question}

By applying the inequality $|V(H)| \leq \alpha(H) \cdot \chi(H)$ to each subgraph $H$ induced by a bag, it is easy to see that $\tw(G)+1 \leq \alpha(G) \cdot \treex(G)$ and $\tw(G)+1 \leq \treea(G) \cdot \chi(G)$.
So if $G$ is a graph with $\alpha(G)=\treea(G)$ or $\chi(G)=\treex(G)$, then $\twG+1 \leq \treeaG \cdot \treexG$ and $G$ is a positive instance of \cref{que:main-question}.
In particular, perfect graphs are positive instances of \cref{que:main-question} since $\chi(G) = \omega(G) \leq \treex(G) \leq \chi(G)$ for every perfect graph $G$.

%\cref{que:main-question} is inspired by the inequality
%\begin{equation}
%    |V| \leq \alpha(G) \cdot \chi(G),\label{eq:simple-alpha-chi}
%\end{equation}
%which holds for every graph $G = (V,E)$.
%In fact, $V$ can be partitioned into $\chi(G)$ independent sets, each of which contains at most $\alpha(G)$ vertices.
%Applying \eqref{eq:simple-alpha-chi} to every bag of a tree decomposition $\T$ of $G$ with $\alpha(\T) = \treeaG$, gives $\twG+1 \leq \treeaG \cdot \chi(G)$.
%Applying \eqref{eq:simple-alpha-chi} to every bag of a tree decomposition $\T$ of $G$ with $\chi(\T) = \treexG$, gives $\twG+1 \leq \alpha(G) \cdot \treexG$.
%Dallard, Milani\v{c}, and \v{S}torgel asked \cref{que:main-question} as a potential common generalization of these two inequalities.
%They observed that the inequality holds for all \emph{perfect graphs}, i.e., for all graphs $G$ where $\chi(H) = \omega(H)$ holds for all induced subgraphs $H$ of $G$.
%In particular, \cref{que:main-question} is confirmed for all chordal graphs and all bipartite graphs.

%In fact, whenever $G$ admits a tree decomposition $\T$ that simultaneously has $\alpha(\T) = \treeaG$ and $\chi(\T) = \treexG$, then applying \eqref{eq:simple-alpha-chi} to each bag shows that $\twG+1 \leq \treeaG \cdot \treexG$ holds for $G$.
%But difficulties arise when the tree decompositions certifying $\treeaG$ and $\treexG$ are distinct and possibly of very different structure.

Even though some positive instances of \cref{que:main-question} are known, we provide a negative answer for \cref{que:main-question} in a strong form in this paper.
(See \cref{theorem:negative-linear-a} for a slightly stronger result.)

\begin{theorem} \label{theorem:negative-linear-a-intro}
    For any functions $f\colon {\mathbb R} \rightarrow {\mathbb R}_{>0}$ and $h\colon {\mathbb R} \rightarrow {\mathbb R}_{>0}$ such that $h(x)$ is non-decreasing and $h(x) = o(x\log\log x)$ as $x \to \infty$,
    there exists a graph $G$ such that
    \[
        \tw(G) > h(\treea(G)) \cdot f(\treex(G))\text{.}
    \]
\end{theorem}

\cref{theorem:negative-linear-a-intro} shows that no upper bound of $\tw(G)$ in terms of $\treea(G)$ and $\treex(G)$ can exist even if we allow the dependence of $\treea(G)$ to grow slightly faster than linear functions and the dependence of $\treex(G)$ to grow arbitrarily fast. %$\treea(G) \log_2\log_2 \treea(G)$.
We complement this result by the following simple observation showing that an upper bound exists if we allow the dependence of $\treea(G)$ to be quadratic.

\begin{proposition} \label{prop:alpha2chi}
  Every graph $G$ satisfies
  $$\tw(G)+1 \leq \treea(G)^{2} \cdot \treex(G)\text{.}$$
\end{proposition}

\begin{proof}
  Take a tree-decomposition whose every bag induces a subgraph with independence number at most $\treea(G)$.
  It suffices to show that every bag contains at most $\treea(G)^{2} \cdot \treex(G)$ vertices.
  For every node $t$ of the tree in the tree-decomposition, let $H_t$ be the subgraph induced by the bag at $t$.
  Note that for every $t$, $\alpha(H_t) \leq \treea(G)$ and $\treex(H_t) \leq \treex(G)$, so $H_t$ has a tree-decomposition whose every bag induces a subgraph of chromatic number at most $\treex(H_t) \leq \treex(G)$ and hence contains at most $\alpha(H_t) \cdot \treex(G) \leq \treea(G) \cdot \treex(G)$ vertices.
  This implies that $\tw(H_t) \leq \treea(G) \cdot \treex(G)-1$.
  It is well known that for every $k$, if $\tw(G) \leq k$, then $\chi(G) \leq k + 1$.
  So $\chi(H_t) \leq \treea(G) \cdot \treex(G)$.
  Hence $|V(H_t)| \leq \alpha(H_t) \cdot \chi(H_t) \leq \treea(G)^2 \cdot \treex(G)$. 
\end{proof}

% \medskip

We remark that in the proof of \cref{prop:alpha2chi}, we showed that $\tw(H_t) \leq \treea(G) \cdot \treex(G)-1$ for every node $t$.
It implies the following correction of \cref{que:main-question} in terms of tree-treewidth.
(Recall that ${\rm tree}\textnormal{-}\tw(G)$ is the minimum $k$ such that there exists a tree-decomposition of $G$ such that every subgraph induced by a bag has treewidth at most~$k$.)

\begin{proposition} \label{prop:treetw-intro}
    For every graph $G$,
    $${\rm tree}\textnormal{-}\tw(G)+1 \leq \treea(G) \cdot \treex(G)\text{.}$$
\end{proposition}

%that $${\rm tree}\textnormal{-}\tw(G)+1 \leq \treea(G) \cdot \treex(G)$$ for every graph $G$, where ${\rm tree}\textnormal{-}\tw(G)$ is the tree-treewidth of $G$.
%\emph{Tree-treewidth} was introduced in \cite{LNW-24}, defined to be the minimum $k$ such that there exists a tree-decomposition such that every bag induces a subgraph of treewidth at most $k$.
Clearly, ${\rm tree}\textnormal{-}\tw(G) \leq \tw(G)$ for every graph $G$.
In general, treewidth cannot be upper bounded by any function of tree-treewidth \cite{LNW-24}.
%But $\tw(G) \leq 150 \cdot {\rm tree}\textnormal{-}\tw(G)$ for almost all $d$-regular graphs for every even integer $d \geq 146$ \cite{LNW-24}.
On the other hand, \cref{prop:treetw-intro} is tight for chordal graphs $G$ since the maximum clique of $G$ must be contained in a bag of every tree-decomposition, implying ${\rm tree}\textnormal{-}\tw(G) +1 \geq \omega(G) = \treex(G) = \treea(G) \cdot \treex(G)$ for every chordal $G$.

This paper is organized as follows.
In \cref{sec:construction}, we introduce a construction that enlarges a graph $G$ to a graph $G'$ without changing the tree-chromatic number while $\tw(G')$ and $\treea(G')$ are closely related to $V(G)$ and $\alpha(G)$.
We then use this construction to provide a negative answer to \cref{que:main-question} in \cref{sec:solution}.

%\subsection{Our results and organization of the paper}
%
%We answer \cref{que:main-question} in the negative by providing a set of counterexamples.
%To this end, we present an infinite family of graphs $G$ with $\treexG = 2$ but $\treeaG \approx \frac14 \twG$.
%The upper bounds on $\treexG$ and $\treeaG$ are certified by two explicit tree decompositions $\T_1$, $\T_2$ of $G$, respectively.
%In $\T_1$, each bag induces a bipartite subgraph of $G$, while $G$ itself has chromatic number $\chi(G) \geq \frac12 \log |V(G)|$.
%However, some bags in $\T_1$ contain large independent sets of $G$.
%On the other hand, every bag in $\T_2$ has small independence number but large chromatic number.
%
%After some preliminaries in \cref{sec:preliminaries}, we introduce so-called $H$-subdivisions of graphs in \cref{sec:construction}.
%We present sharp bounds on the treewidth, tree-independence number, and tree-chromatic number of $H$-subdivisions of complete graphs.
%Taking suitable graphs for $H$, we obtain our counterexamples to \cref{que:main-question}; see \cref{thm:random-construction} for the best lower bounds.
%
%In \cref{sec:upper-bounds} we give several upper bounds on $\twG$ in terms of $\treeaG$ and $\treexG$ by employing results from Ramsey theory.
%There we also discuss limitations of our approach from \cref{sec:construction} and pose some weakenings of \cref{que:main-question}.
%We conclude with more open problems in \cref{sec:conclusions}.

%\section{Preliminaries}
%\label{sec:preliminaries}
%

\section{Completion of graphs}
\label{sec:construction}

In this section we provide a key ingredient for constructing a negative answer to \cref{que:main-question}. 
%In this section we break down in detail the construction that allows us to prove our main theorem.
%More precisely, by constructing an infinite family $\mathcal G$ of graphs such that each $G\in \mathcal{G}$ satisfies the inequality $\tw(G)+1>\treeaG \cdot \treexG$, we answer in the negative \cref{que:main-question}.
%We begin by formally stating our main result.

%\begin{theorem}\label{thm-deterministic-construction}
%    There exists an infinite sequence of graphs $(G_k)_{k \in \mathbb{N}_+}$ such that:
%    \begin{enumerate}
%        \item Each graph $G_k$ satisfies the inequality
%        \[{\rm tw}(G_k)+1>\treea(G_k) \cdot \treex(G_k).\]
%        \item The ratio between $\tw(G_k)$ and $\treea(G_k) \cdot \treex(G_k)$ converges to a number that is at least $\frac{4}{3}$, that is:
%        \[\lim_{k\to \infty}\frac{\tw(G_k)}{\treea(G_k) \cdot \treex(G_k)} \geq \frac{4}{3}.\]
%        \item Each graph $G_k$ in the sequence can be constructed deterministically in linear time. \label{thm-item:deterministic-construction-3}
%    \end{enumerate}
%\end{theorem}

For a graph $H$, the {\it 1-completion} of $H$ is the graph obtained from $H$ by, for each pair of non-adjacent vertices $u$ and $v$ of $H$, adding a new vertex of degree~$2$ adjacent to exactly $u$ and $v$.
In this paper, denote the 1-completion of $H$ by $C(H)$.

The first observation is that taking the 1-completion does not change the tree-chromatic number.
% It is well known that the parameters $\tw, \treea, \treex$ are monotone under taking minors, induced minors, and subgraphs, respectively.
% For convenience, we state this as a separate lemma.

% \begin{lemma}\label{lemma:relationship-monotonicities}
%    Let $G$ be a graph. Then the following inequalities hold.
%    \begin{itemize}
%        \item For any minor $H$ of $G$, $\tw(H)\leq \tw(G)$ (see \cite{RSP-86}). 
%        \item For any induced minor $H$ of $G$, $\treea(H)\leq \treea(G)$ (see \cite{DMS-24}).
%        \item For any subgraph $H$ of $G$, $\treex(H)\leq \treex(G)$ (see \cite{Sey-16}). 
%    \end{itemize}
% \end{lemma}

\begin{lemma}\label{lemma:completion-treex}
    If $H$ is a graph with $E(H) \neq \emptyset$, then $\treex(C(H)) = \treex(H)$.
\end{lemma}
\begin{proof}
    Since $H$ is a subgraph of $C(H)$, we get that $\treex(C(H))\geq \treex(H)$. %by \cref{lemma:relationship-monotonicities}.
    %($\treex$ is monotone under taking subgraphs\cite{Sey-16}.)
    
    To prove $\treex(C(H))\leq \treex(H)$, it suffices to construct a tree-decomposition $\T$ of $C(H)$ with $\chi(\T) = \treex(H)$.

    For any pair of non-adjacent vertices $u$ and $v$ in $H$, let $a_{uv}$ be the vertex of $C(H)-V(H)$ adjacent to both $u$ and $v$.
    Let $\T_H = (T_H, \{X_t\}_{t\in V(T_H)})$ be a tree-decomposition of $H$ with $\chi(\T_H) = \treex(H)$. 
    % Note that for every edge $uv$ of $H$, there exists $t_{uv} \in V(T_H)$ such that $X_{t_{uv}} \supseteq \{u,v\}$.
    
    % Define $T$ to be the tree obtained from $T_H$ by, for each edge $uv$ of $H$, adding a leaf $t'_{uv}$ adjacent to $t_{uv}$.
    % For every $t \in V(T)-V(T_H)$, we know $t=t_{uv}'$ for some $uv \in E(H)$, and we let $Y_{t'_{uv}}=\{u,v,a_{uv}\}$.
    % For every $t \in V(T_H)$, let $Y_t = X_t \cup \{a_{uv}: u,v$ are distinct vertices of $H$ with $uv \not \in E(H)\}$.
    % Let $\T = (T, (Y_t)_{t \in V(T)})$.
    % Clearly, $\T$ is a tree-decomposition of $C(H)$.

    We construct a new tree-decomposition from $\T_H$.
    Fix $u, v \in V(H)$ such that $uv \notin E(H)$.
    Either (1) there exists $t_{uv} \in V(T_H)$ such that $\{u,v\} \subseteq X_{t_{uv}}$; or (2) there is no bag in $\T_H$ containing both $u$ and $v$.
    In case (1), we add a new vertex $t'_{uv}$ adjacent to $t_{uv}$, and let $X_{t'_{uv}}=\{u,v,a_{uv}\}$.
    In case (2), we add $a_{uv}$ into $X_t$ for every $t \in V(T_H)$; note that $a_{uv}$ has degree at most 1 in the subgraph induced by the new $X_t$. 
    Clearly, this procedure yields a tree-decomposition $\T = (T, \{Y_t\}_{t \in V(T)})$ of $C(H)$.

    Since $E(H) \neq \emptyset$, some bag of $\T_H$ contains an edge, so $\chi(\T_H) \geq 2$.
    For every $t \in V(T)-V(T_H)$, the subgraph of $C(H)$ induced by $Y_t$ is a forest and hence has chromatic number at most $2 \leq \chi(\T_H)$.
    For every $t \in V(T_H)$, the subgraph of $C(H)$ induced by $Y_t$ is obtained from $H[X_t]$ by adding isolated vertices or attaching leaves and hence has chromatic number at most $\max\{2,\chi(H[X_t])\} \leq \chi(\T_H)$.
    So $\chi(\T) \leq \chi(\T_H) = \treex(H)$.
\end{proof}

For every graph $G$, we denote the complement of $G$ by $\overline{G}$.

The second observation is that the treewidth and the tree-independence number of the 1-completion of $H$ is closely related to $|V(H)|$ and the independence number of $H$, respectively.

\begin{lemma}\label{lemma:completion-tw}
    If $H$ is a graph with $|V(H)| \geq 3$, then $\tw(C(H)) = |V(H)|-1$ and $\treea(C(H)) \leq \alpha(H)$.
\end{lemma}
\begin{proof}
    Let $n=|V(H)|$.
    Since $C(H)$ contains $K_n$ as a minor, $\tw(C(H)) \geq \tw(K_n) = n-1$. %by \cref{lemma:relationship-monotonicities}.
    
    To prove $\tw(C(H)) \leq n-1$, we construct a tree decomposition $\T = (T,\{X_t\}_{t\in V(T)})$ of $C(H)$ of width $n-1$. 
    Let $T=K_{1,|E(\overline H)|}$.
    Let the bag of the central node of $T$ be $V(H)$. 
    Note that there exists a bijection $\sigma$ from the set of leaves of $T$ to $E(\overline{H})$.
    For every leaf $t$ of $T$, let $X_t$ be consisting of the ends of $\sigma(t)$ and the vertex of $C(H)-V(H)$ adjacent to both ends of $\sigma(t)$.
    Then $\T = (T,\{X_t\}_{t\in V(T)})$ is a tree-decomposition of $C(H)$ of width $\max\{|V(H)|,3\}-1 = n-1$. 

    Moreover, for every $t \in V(T)$, the subgraph of $C(H)$ induced by $X_t$ is either $H$ or a path on three vertices, so $\alpha(\T) = \max\{\alpha(H),2\}$.
    If $\alpha(H) \geq 2$, then $\alpha(\T) = \alpha(H)$.
    If $\alpha(H) \leq 1$, then $H$ is a complete graph and $T$ consists of a single vertex, so $\alpha(\T)=\alpha(H)$.
    Hence $\treea(C(H)) \leq \alpha(\T) = \alpha(H)$.
%    Note that the induced subgraph of $S(\overline H)$ on vertex set $V$ is exactly $H$, since $S(\overline H)$ arises from $K_n$ by subdividing the edges of $\overline{H}$, i.e., non-edges of $H$.
%    Now, label each of the pendant nodes of $T$ by an edge $e\in E(\overline{H})$ bijectively and let the bag corresponding to the pendant node in $T$ labeled $e=uv$ consist of vertices $\{u,v,e\}\subseteq %V(S(\overline H))$.
%    Observe that this construction yields a valid tree decomposition $\T$ of $S(\overline H)$ of width $n-1$, proving \cref{lemma-item:treewidth}. 
\end{proof}

Even though \cref{lemma:completion-tw,lemma:completion-treex} are sufficient to construct a simple negative instance of \cref{que:main-question}\footnote{E.g. by taking $H$ to be isomorphic to $k$ copies of $C_5$ (for any constant $k\ge 1$), then by \cref{lemma:completion-treex}, we have $\treex(C(H)) = 2$ and by \cref{lemma:completion-tw}, it follows that $\tw(C(H)) = 5k-1$ and $\treea(C(H))\le 2k$.}, in the rest of this section we strengthen \cref{lemma:completion-tw} to show that the upper bound for $\treea(C(H))$ is almost tight.

%Before proving \cref{thm-deterministic-construction}, we first introduce a useful auxiliary construction of graphs, and prove some properties that such graphs satisfy.
%This will allow us to easily construct a class that satisfies the conditions of \cref{thm-deterministic-construction}.

%Let $G = (V,E_G)$ be a graph and $H = (V,E_H)$ be a spanning subgraph of $G$.
%Then we define the \emph{$H$-subdivision of $G$}, denoted by $S_G(H)$, as the graph obtained from $G$ by subdividing each edge $e$ that is present in $H$.
%Formally
%\begin{align*}
%    V(S_G(H)) &= V \cup E_H\\
%    E(S_G(H)) &= \left(E_G \setminus E_H \right) \cup \{\{u,e\}\mid e\in E_H, u\in e \}
%\end{align*}
%For the purpose of our construction, given a graph $H=(V,E_H)$ the canonical choice for the graph $G$ will be the graph isomorphic to the complete graph $K_n$ on vertex set $V$.
%Thus, if $G$ is a complete graph, we omit the index $G$ in the notation for the $H$-subdivision of $G$, and simply write $S(H)$.
%We will soon demonstrate how we can use $H$-subdivisions of complete graphs to prove our main theorem.
%But before that, let us prove the following lemma.

\begin{lemma}\label{lemma:subdivision-of-Kn-tree-alpha-lower-bound}
    For every $n \geq 2$ we have $\treea(C(\overline{K_n})) = n-1$. %$\treea(S(K_n)) = n-1$.
\end{lemma}
\begin{proof}
    Let $A \subset V(C(\overline{K_n}))$ be the subset of $n$ vertices in $C(\overline{K_n})$ corresponding to the original vertices of $\overline {K_n}$.
    Note that $A$ is an independent set in $C(\overline{K_n})$.
    
    We first show the lower bound $\treea(C(\overline{K_n})) \geq n-1$.
    Consider a tree decomposition $\T = (T,\{X_t\}_{t\in V(T)})$ of $C(\overline{K_n})$ with $\alpha(\T) = \treea(C(\overline{K_n}))$.
    We may assume that no bag of $\T$ contains $A$, for otherwise $\alpha(\T) \geq |A| = n \geq n-1$, as desired.
    %If for some node $t \in V(T)$ it holds $A \subseteq X_t$, then $\alpha(\T) \geq |A| = n \geq n-1$, as desired.
    So there exists an edge $t_1t_2$ of $T$ such that $T - t_1t_2$ consists of two subtrees $T_1$ and $T_2$ of $T$, and there exists a partition $A_1 \cup A_2 \cup A_3$ of $A$ with the following properties:
    \begin{itemize}
        \item $|A_1| \geq 1$ and every vertex in $A_1$ appears only in bags of $T_1$.
        \item $|A_2| \geq 1$ and every vertex in $A_2$ appears only in bags of $T_2$.
        \item Every vertex in $A_3$ appears in bags of $T_1$ and bags of $T_2$.
    \end{itemize}
    Clearly, $|A_1|+|A_2|+|A_3| = |A| = n$.
    Moreover, the bag $X_{t_1}$ contains all vertices in $A_3$, as well as the set $B$ consisting of the (degree-$2$) vertices in $C(\overline{K_n})-A$ adjacent in $C(\overline{K_n})$ to both vertices in $A_1$ and $A_2$.
    Observe that $A_3 \cup B$ is an independent set in $C(\overline{K_n})$.
    So $\alpha(\T) \geq |A_3 \cup B| = |A_3| + |B|$.
    Since $|B| = |A_1|\cdot |A_2| \geq |A_1|+|A_2|-1$, where the inequalities uses that $|A_1|,|A_2| \geq 1$, we have $\alpha(\T) \geq |B| + |A_3| \geq |A_1|+|A_2|+|A_3|-1 = n-1$, as desired.

    \medskip

    Now we show the upper bound $\treea(C(\overline{K_n})) \leq n-1$.
    The case $n=2$ is trivial, so we may assume $n \geq 3$.
    Let $S$ be the star with $\binom{n-1}{2}+1$ leaves.
    We denote by $s$ the central vertex of $S$, and denote the leaves of $S$ by $s_0$ and $s_{i,j}$ for $1 \leq i <j \leq n-1$.
    Denote $A$ by $\{a_0,a_1,a_2,...,a_{n-1}\}$.
    Define $Y_s = (A-\{a_0\}) \cup N_{C(\overline{K_n})}(a_0)$ and $Y_{s_0} = \{a_0\} \cup N_{C(\overline{K_n})}(a_0)$.
    For every $i,j$ with $1 \leq i <j \leq n-1$, define $Y_{s_{i,j}}$ to be the set consisting of $a_i,a_j$ and the unique common neighbor of $a_i$ and $a_j$ in $C(\overline{K_n})$.
    Clearly, $(S,(Y_t)_{t \in V(S)})$ is a tree-decomposition of $C(\overline{K_n})$.
    Note that the subgraph induced by $Y_s$ contains a perfect matching with $n-1$ edges; the subgraph induced by each leaf of $S$ is a star with at most $\max\{|A|-1,2\} = n-1$ leaves.
    So every bag induces a subgraph with independence number at most $n-1$.
    This shows $\treea(C(\overline{K_n})) \leq n-1$.
%    , we construct a tree decomposition $\T = (T,\{X_t\}_{t\in V_T})$ of $S(K_n)$ with $\alpha(\T) = n-1$.
%    For this, pick one vertex $a \in A$ and let $F$ be the set of the $n-1$ edges in $K_n$ incident to $a$.
%    For $T$ we take a star with $\binom{n-1}{2}+1$ leaves and denote by $t$ the central vertex of $T$.
%    For the corresponding bag $X_t$, we take the union of $A - \{a\}$ and all subdivision vertices corresponding to edges in $F$.
%    Note that $S(K_n)[X_t]$ is a matching on $n-1$ edges.
%    
%    For each edge $e = uv$ in $E(K_n) - F$, let $t(e)$ be a corresponding leaf in $T$.
%    We define the corresponding bag as $X_{t(e)} = \{u,v,e\}$, where $e$ is the subdivision vertex of the edge $uv$ in $S(K_n)$.
%    Note that $S(K_n)[X_{t(e)}]$ is a path on two edges.
%    
%    Finally, we incorporate the special vertex $a \in A$ and its incident edges.
%    For this, let $t(a)$ be the last leaf of $T$ and define $X_{t(a)} = \{a\} \cup F$.
%    That is, $X_{t(a)}$ contains the original vertex $a$ of $K_n$ and the subdivision vertex of each edge $e$ incident to $a$ in $K_n$.
%    Note that $S(K_n)[X_{t(a)}]$ is a star with $n-1$ leaves.
%
%    Finally, observe that $\T = (T,\{X_t\}_{t\in V_T})$ is a valid tree decomposition of $S(K_n)$, and $\alpha(\T) = n-1$.
\end{proof}

%We can now prove some nice properties that the graph $S(H)$ satisfies.
%As it will be convenient later, let us actually phrase these properties in terms of the complement graph.

Combining \cref{lemma:completion-tw,lemma:subdivision-of-Kn-tree-alpha-lower-bound}, we can determine $\treea(C(H))$ quite precisely.

\begin{lemma}\label{lemma:completion-treea}
    If $H$ is a graph with $|V(H)| \geq 3$, then $\alpha(H)-1 \leq \treea(C(H)) \leq \alpha(H)$.
\end{lemma}
\begin{proof}
    By \cref{lemma:completion-tw}, $\treea(C(H)) \leq \alpha(H)$.
    Since $C(H)$ contains $C(\overline{K_{\alpha(H)}})$ as an induced subgraph, the inequality $\treea(C(H)) \geq \alpha(H)-1$ follows directly from \cref{lemma:subdivision-of-Kn-tree-alpha-lower-bound}, whenever $\alpha(H) \geq 2$.
    On the other hand, if $\alpha(H) \leq 1$, we have $\treea(C(H)) \geq 0 \geq \alpha(H)-1$ trivially.
\end{proof}

\section{Negative answer to \cref{que:main-question}} \label{sec:solution}

We shall use a version of \emph{Burling graphs} \cite{burling1965coloring} to construct a negative answer to \cref{que:main-question}.
This is a family of graphs, which are recursively defined as follows.
Throughout this section, $(G_{n}, \mathcal{S}_{n})_{n \in {\mathbb N}}$ is a sequence of ordered pairs, where $G_{n}$ is a graph and $\mathcal{S}_{n}$ is a collection of stable sets in $G_n$ for every $n \in {\mathbb N}$, such that $G_{1} = K_{1}$ and $\mathcal{S}_{1} = \{V(G_{1})\}$, and for $n \geq 2$, $G_n$ and $\mathcal{S}_n$ are defined as follows:
\begin{itemize}
  %\item Create a copy \(G\) of \(G_{n-1}\).
  \item For each $S \in \mathcal{S}_{n-1}$, create a copy $G_{S}$ of $G_{n-1}$, and let $\mathcal{S}_S$ be the collection of stable sets of $G_S$ consisting of the copies of the members of $\mathcal{S}_{n-1}$ in $G_S$.
  \item Define $G_n$ to be the graph obtained from $G_{n-1} \cup \bigcup_{S \in \mathcal{S}_{n-1}}G_S$, by for each $S \in \mathcal{S}_{n-1}$ and $Q \in \mathcal{S}_{S}$, adding a vertex $v_{S,Q}$ adjacent to all vertices in $Q$.
  \item Define $\mathcal{S}_{n} = \{S \cup Q, S \cup \{v_{S,Q}\}: S \in \mathcal{S}_{n-1}, Q \in \mathcal{S}_S\}$. 
\end{itemize}
%A graph is called a \emph{Burling graph} if it is an induced subgraph of a graph in the Burling sequence.

First, we prove that for each graph $G_n$ in the above sequence we have $\treex(G_n) \leq 2$.

A \emph{star forest} is a graph whose every component is a star.
Note that $K_1$ is a star.

\begin{lemma} \label{lemma:Burling-treex}
For every positive integer $n$, there exists a tree-decomposition \(\mathcal{T}\) of $G_n$ such that every bag of $\T$ induces a star forest, and every member of $\mathcal{S}_n$ is contained in a bag of $\T$. %for every $S \in \mathcal{S}_n$, some bag of $\mathcal{T}$ contains $S$.
\end{lemma}

\begin{proof}
We shall prove this lemma by induction on $n$.
The case $n=1$ is obvious.
So we may assume $n \geq 2$.
By the inductive hypothesis, there exists a tree-decomposition $\T_{n-1}=(T_{n-1},\{X_t\}_{t \in V(T_{n-1})})$ of $G_{n-1}$ such that every bag of $\T_{n-1}$ is a star forest, and for every $S \in \mathcal{S}_{n-1}$, there exists $t_S \in V(T_{n-1})$ with $S \subseteq X_{t_S}$.

For every $S \in \mathcal{S}_{n-1}$, let $\T_S = (T_S, \{X^S_t\}_{t \in V(T_S)})$ be a tree-decomposition $G_{S}$ isomorphic to $\T_{n-1}$. %, and let $\T_S'$ be the tree-decomposition $(T_S, (X^S_t)_{t \in V(T_S)})$ of $G_S \cup G_{n-1}[S]$ such that obtained from $\T_S$ by adding $S$ into every bag of $\T_S$.
Note that for any $S \in {\mathcal S}_{n-1}$ and $Q \in \mathcal{S}_S$, there exists $t_{S,Q} \in V(T_S)$ such that $Q \subseteq X^S_{t_{S,Q}}$.
For every $S \in {\mathcal S}_{n-1}$, we construct a tree $T_S'$ by for every $Q \in {\mathcal S}_S$, adding a new node $t'_{S,Q}$ into $T_S$ and the edge $t'_{S,Q}t_{S,Q}$, and define $Y_{t'_{S,Q}}=Q \cup \{v_{S,Q}\}$.
For every $S \in {\mathcal S}_{n-1}$, let $Y_{t} = X^S_t \cup S$ for every $t \in V(T_S)$.

%Note that for every $S \in {\mathcal S}_{n-1}$, there exists $t_S \in V(T_S)$ such that $S \subseteq X^S_{t_S}$.
Let $T_n$ be the tree obtained from $T_{n-1} \cup \bigcup_{S \in \mathcal{S}_{n-1}}T_S'$ by, for each $S \in {\mathcal S}_{n-1}$, adding an edge between $t_S$ and some node of $T_S'$.
Then $(T_n, \{X_t\}_{t \in V(T_{n-1})} \cup \{Y_t\}_{S \in {\mathcal S}_{n-1}, t \in V(T_S')})$ is a tree-decomposition of $G_n$ such that every member of $\mathcal{S}_n$ is contained in some bag.

Moreover, for every $t \in V(T_{n-1})$, we know $G_n[X_t]=G_{n-1}[X_t]$ is a star forest by the inductive hypothesis; for any $S \in {\mathcal S}_{n-1}$ and $t \in V(T_S')$, if $t=t'_{S,Q}$ for some $Q \in \mathcal{S}_S$, then $G_n[Y_t]$ is a star, and if $t \in V(T_S)$, then $Y_t$ is a union of a subset $X^S_t$ of $V(G_S)$ and a subset $S$ of $V(G_{n-1})$, so $Y_t$ is a stable set in $G_n$.
Therefore, every bag of $(T_n, \{X_t\}_{t \in V(T_{n-1})} \cup \{Y_t\}_{S \in {\mathcal S}_{n-1}, t \in V(T_S')})$ induces a star forest in $G_n$.
\end{proof}

Walczak \cite{Wal-15} proved that Burling graphs have unbounded fractional chromatic number.
This result can be stated in terms of its dual.
For a graph $G$, a function $f\colon V(G) \rightarrow {\mathbb R}$, and a subset $S$ of $V(G)$, we define $f(S)$ to be $\sum_{x \in S}f(x)$.

\begin{lemma}[Walczak {\cite[Page 222]{Wal-15}}] \label{lemma:Burling-frac}
    For every positive integer $k$, there exists a function $w_k\colon V(G_k) \rightarrow {\mathbb N}$ such that $w_k(V(G_k)) = \frac{k+1}{2} \cdot 2^{2^{k-1}-1}$, and  $w_k(I) \leq 2^{2^{k-1}-1}$ for every stable set $I$ of $G_k$.
\end{lemma}

For two graphs $G$ and $H$ a map $f \colon V(G) \to V(H)$ is a \emph{(graph) homomorphism} from $G$ to $H$ if $f(u)f(v) \in E(H)$ for every $uv \in E(G)$.
It is well-known (see for example \cite[Corollary 1.8]{HN-04}) that if there exists a homomorphism from $G$ to $H$, then $\chi(G) \leq \chi(H)$.
An analogous statement holds for the tree-chromatic number.

\begin{lemma} \label{lemma:homo-treex}
Let $G$ and $H$ be graphs.
If there exists a homomorphism from $G$ to $H$, then $\treex(G) \leq \treex(H)$.
\end{lemma}

\begin{proof}
Let $\T=(T,\{X_t\}_{t \in V(T)})$ be a tree-decomposition of $H$ such that $\chi(\T)=\treex(H)$.
Let $f$ be a homomorphism from $G$ to $H$.
For every $t \in V(T)$, let $X'_t = \bigcup_{v \in X_t}f^{-1}(\{v\})$.
Note that $(T,\{X'_t\}_{t \in V(T)})$ is a tree-decomposition of $G$.

Since $f$ is a homomorphism, $\bigcup_{x \in S}f^{-1}(\{x\})$ is a stable set in $G$ for every stable set $S$ in $H$.
Hence $\chi(G[X'_t]) \leq \chi(H[X_t]) \leq \chi(\T) = \treex(H)$ for every $t \in V(T)$.
\end{proof}

Next, we construct a blowup of each $G_k$ based on function $w_k$ in \cref{lemma:Burling-frac}.
That is, each vertex $v \in V(G_k)$ is blown up to a set of $w_k(v)$ pairwise non-adjacent vertices.
We summarize the crucial properties of the resulting graph $H_k$.

\begin{lemma} \label{lemma:blowup-Burling}
For every positive integer $k$, there exists a graph $H_k$ such that 
    \begin{enumerate}
        \item $|V(H_k)| = \frac{k+1}{2} \cdot 2^{2^{k-1}-1}$,
        \item $\alpha(H_k) \leq \frac{2}{k+1}|V(H_k)| < \frac{4|V(H_k)|}{\log_2\log_2|V(H_k)|}$, and 
        \item $\treex(H_k) = 2$.
    \end{enumerate}
\end{lemma}
\begin{proof}
Let $k$ be a positive integer.
By \cref{lemma:Burling-frac}, there exists a function $w_k\colon V(G_k) \rightarrow {\mathbb N}$ such that $w_k(V(G_k)) = \frac{k+1}{2} \cdot 2^{2^{k-1}-1}$, and $w_k(I) \leq 2^{2^{k-1}-1}$ for every stable set $I$ of $G_k$.

For every $v \in V(G_k)$, let $S_v$ be a set consisting of $w_k(v)$ pairwise non-adjacent vertices.
Define $H_k$ to be the graph obtained from the disjoint union of $S_v$ (for all $v \in V(G_k)$) by adding edges such that every vertex in $S_a$ is adjacent in $H_k$ to every vertex in $S_b$ for any $ab \in E(G_k)$.
%Define $H_k$ to be the graph obtained from $G_k$ by, for every $v \in V(G_k)$, replacing $v$ by a stable set $S_k$ of size $w_k(v)$, such that every vertex in $S_a$ is adjacent in $H_k$ to every vertex in $S_b$ for any $ab \in E(G_k)$.
Then $|V(H_k)| = w_k(V(G_k)) = \frac{k+1}{2} \cdot 2^{2^{k-1}-1}$.

Let $A$ be a stable set in $H_k$.
By the definition of $H_k$, there exists no $ab \in E(G_k)$ such that $A \cap S_a \neq \emptyset$ and $A \cap S_b \neq \emptyset$.
Hence there exists a stable set $I_A$ of $G_k$ such that $A \subseteq \bigcup_{v \in I_A}S_v$.
So $|A| \leq \sum_{v \in I_A}|S_v| = w_k(I_A) \leq 2^{2^{k-1}-1}$.

Therefore, $\alpha(H_k) \leq 2^{2^{k-1}-1} = \frac{2}{k+1}|V(H_k)|$.
Since $k \geq 1$, we have $\log_2\log_2(k+1) \leq k$.
Hence $$\log_2\log_2|V(H_k)| = \log_2(\log_2(k+1)-1) + \log_2(2^{k-1}-1) < \log_2\log_2(k+1) + k-1 \leq 2k.$$
Therefore, $\alpha(H_k) \leq \frac{2}{k+1}|V(H_k)| < \frac{4|V(H_k)|}{\log_2\log_2|V(H_k)|}$.

By \cref{lemma:Burling-treex}, there exists a tree-decomposition $\T$ of $G_k$ such that every bag induces a star forest.
So $\treex(G_k) \leq \chi(\T) \leq 2$.
Since the function $\phi \colon V(H_k) \rightarrow V(G_k)$ mapping each vertex $h \in V(H_k)$ to the vertex $v \in V(G_k)$ such that $h \in S_v$ is a homomorphism from $H_k$ to $G_k$, we have $\treex(H_k) \leq \treex(G_k) \leq 2$ by \cref{lemma:homo-treex}.
Since $E(H_k) \neq \emptyset$, every tree-decomposition of $H_k$ has a bag that contains an edge, so $\treex(H_k) \geq 2$.
This shows $\treex(H_k) = 2$.
\end{proof}

Using the 1-completion $C(H_k)$ of the graphs $H_k$ from \cref{lemma:blowup-Burling}, we are ready to prove our main result from which \cref{theorem:negative-linear-a-intro} will easily follow.

\begin{theorem} \label{theorem:negative-linear-a}
    For any functions $f\colon {\mathbb R} \rightarrow {\mathbb R}_{>0}$ and $g\colon {\mathbb R} \rightarrow {\mathbb R}_{>0}$ with $g(x) = o(\log\log x)$ as $x \to \infty$, there exists a graph $G$ such that $$\tw(G) > \treea(G) \cdot f(\treex(G)) \cdot g(\tw(G)).$$
\end{theorem}

\begin{proof}
Suppose to the contrary that $\tw(G) \leq \treea(G) \cdot f(\treex(G)) \cdot g(\tw(G))$ for every graph $G$.
%Since $\tw(K_2)>0$ and $\treex(K_2)=2$, we know $f(2) > 0$.

For every positive integer $k$, let $H_k$ be the graph mentioned in \cref{lemma:blowup-Burling}.
Then $\tw(C(H_k)) \leq \treea(C(H_k)) \cdot f(\treex(C(H_k))) \cdot g(\tw(C(H_k)))$. 
When $k \geq 2$, we have $|V(H_k)| \geq 3$, so \cref{lemma:completion-treex,lemma:completion-tw} imply
    \begin{align*}
        |V(H_k)|-1 = \tw(C(H_k)) & \leq \treea(C(H_k)) \cdot f(\treex(C(H_k))) \cdot g(\tw(C(H_k))) \\
        & \leq \alpha(H_k) \cdot f(\treex(H_k)) \cdot g(|V(H_k)|-1).
    \end{align*}
By \cref{lemma:blowup-Burling}, 
    %\begin{align*}
        $$|V(H_k)|-1 \leq \alpha(H_k) \cdot f(\treex(H_k)) \cdot g(|V(H_k)|-1) \leq \frac{4|V(H_k)|}{\log_2\log_2|V(H_k)|} \cdot f(2) \cdot g(|V(H_k)|-1).$$
    %\end{align*}
%Hence $\frac{1}{4f(2)}(1 - \frac{1}{|V(H_k)|}) \leq \frac{g(|V(H_k)|-1)}{\log_2\log_2|V(H_k)|} \leq \frac{g(|V(H_k)|-1)}{\log_2\log_2(|V(H_k)|-1)}$.
Since $|V(H_k)| \to \infty$ as $k \to \infty$ by \cref{lemma:blowup-Burling}, we know %$\frac{1}{4f(2)} = \lim_{k \to \infty}\frac{1}{4f(2)}(1 - \frac{1}{|V(H_k)|}) \leq \lim_{k \to \infty}\frac{g(|V(H_k)|-1)}{\log_2\log_2(|V(H_k)|-1)} = 0$, a contradiction.
    \begin{align*}
        0 < \frac{1}{4f(2)} = \lim_{k \to \infty}\frac{1}{4f(2)}(1 - \frac{1}{|V(H_k)|}) & \leq \lim_{k \to \infty}\frac{g(|V(H_k)|-1)}{\log_2\log_2|V(H_k)|} \\
        & \leq \lim_{k \to \infty}\frac{g(|V(H_k)|-1)}{\log_2\log_2(|V(H_k)|-1)} = 0,
    \end{align*}
a contradiction.
\end{proof}

Now we deduce \cref{theorem:negative-linear-a-intro} from \cref{theorem:negative-linear-a}.

\begin{corollary}
    For any functions $f \colon {\mathbb R} \rightarrow {\mathbb R}_{>0}$ and $h \colon {\mathbb R} \rightarrow {\mathbb R}_{>0}$ such that $h(x)/x$ is non-decreasing and $h(x) = o(x\log\log x)$ as $x \to \infty$, there exists a graph $G$ such that $$\tw(G) > h(\treea(G)) \cdot f(\treex(G)).$$
\end{corollary}

\begin{proof}
Let $g(x) = h(x)/x$ for every $x>0$.
Then $g(x) = o(\log\log x)$ as $x \to \infty$.
By \cref{theorem:negative-linear-a}, there exists a graph $G$ such that $\tw(G) > \treea(G) \cdot f(\treex(G)) \cdot g(\tw(G)) \geq \treea(G) \cdot f(\treex(G)) \cdot g(\treea(G)) = h(\treea(G)) \cdot f(\treex(G))$.
\end{proof}

%\input{construction}
%\input{upper-bounds}

%\section{Conclusions}
%\label{sec:conclusions}

    % We compare the treewidth $\twG$ of a graph $G$ to its related parameters $\treeaG$ and $\treexG$.
% In particular,
%We seek to find graphs $G$ with large $\twG$ but small $\treeaG$ and small $\treexG$.
%In order to improve our exponential upper bounds, it might be worthwhile to find the exact asymptotics for constant $\treex$ (or constant $\treea$).

%\begin{question}
%    Do graphs $G$ with $\treexG \leq k$ satisfy $\twG \leq \mathcal{O}_k(\treeaG)$?
%\end{question}

%Our best constructions are $\overline H$-subdivisions of complete graphs, for which we seek graphs $H$ with small $\treex(H)$ and small $\alpha(H)$, but many vertices.
%Having focused on $\treex(H) = 2$, our best examples have $\alpha(H) \approx \frac14|V(H)|$.

%\begin{question}
%    Does every graph $H$ with $\treex(H) = 2$ satisfy $\alpha(H) \geq \frac14 |V(H)|$?
%    % Is it true that $\alpha(H) \geq \frac14 |V(H)|$ for all graphs $H$ with $\treex(H) = 2$?

%    Are there graphs $H$ with $\alpha(H) < \frac{1}{2 \cdot \treex(H)}|V(H)|$?
%\end{question}

\bibliographystyle{plainurl}
\bibliography{references}

\end{document}